\documentclass[11pt]{amsart}
\usepackage{amsmath,amssymb,mathrsfs}
\newtheorem{theorem}{Theorem}
\newtheorem{proposition}[theorem]{Proposition}
\newtheorem{corollary}[theorem]{Corollary}
\newtheorem{lemma}[theorem]{Lemma}
\begin{document}

\title[Spheres of constant mean curvature]{Large outlying stable constant mean curvature spheres in initial data sets}
\author{Simon Brendle and Michael Eichmair}
\begin{abstract}
We give examples of asymptotically flat three-manifolds $(M,g)$ which admit arbitrarily large constant mean curvature spheres that are far away from the center of the manifold. This resolves a question raised by G. Huisken and S.-T. Yau in 1996. On the other hand, we show that such surfaces cannot exist when $(M,g)$ has nonnegative scalar curvature. This result depends on an intricate relationship between the scalar curvature of the initial data set and the isoperimetric ratio of large stable constant mean curvature surfaces.
\end{abstract}
\address{Department of Mathematics \\ Stanford University \\ Stanford, CA 94305 \\ U.S.A.}
\address{Departement Mathematik \\ ETH Z\"urich \\ 8092 Z\"urich \\ Switzerland}
\maketitle

\section{Introduction}

In this paper, we contribute towards a full understanding of large stable constant mean curvature surfaces in asymptotically flat three-manifolds. In particular, we discover and use a delicate mechanism that brings out the subtle role of scalar curvature in this problem. Our results lie at the juncture of the classical analysis of Riemannian manifolds via their isoperimetric properties on the one hand, and deep recent developments in mathematical relativity through the work of H.~Bray \cite{Bray-thesis} and G.~Huisken \cite{Huisken} on the other hand.

Stable constant mean curvature surfaces have been used with great success to capture physical properties of initial data sets for the Einstein equations.
A landmark result in this direction due to G.~Huisken and S.-T.~Yau \cite{Huisken-Yau} says that the ends of non-flat initial data sets admit canonical foliations through stable constant mean curvature spheres. The mass and the center of mass of the initial data set are encoded in these spheres. Moreover, any stable constant mean curvature sphere which encloses a sufficiently large compact set is a leaf of this foliation (cf. \cite{Qing-Tian}). D.~Christodoulou and S.-T.~Yau \cite{Christodoulou-Yau} observed that the Hawking mass is a natural quasi-local measure of gravitation for regions in initial data sets that are bounded by stable constant mean curvature spheres. H.~Bray \cite{Bray-thesis} observed that the Hawking mass is monotone increasing along an expanding foliation through such spheres, such as that of \cite{Huisken-Yau}, if the scalar curvature is non-negative.

J.~Metzger and the second-named author showed that the stable constant mean curvature spheres in \cite{Huisken-Yau} enclose their volume with the least amount of area possible. In fact, they are uniquely characterized by this property  \cite{Eichmair-Metzger1}. This result holds in all dimensions \cite{Eichmair-Metzger3}. It resolves a conjecture of H. Bray who in this Stanford thesis characterized the solutions of the isoperimetric problem in the exact Schwarzschild manifold \cite{Bray-thesis}. We characterized the solutions of the isoperimetric problem in the doubled Schwarzschild manifold in \cite{Brendle-Eichmair}. Moreover, in a recent paper \cite{Brendle2}, the first-named author showed that any embedded constant mean curvature surface in the Schwarzschild manifold is a sphere of symmetry. This result is similar to the classical Alexandrov theorem in Euclidean space in that it does not require any assumptions concerning the topology of the surface or the stability operator. Moreover, the result can be generalized to a large class of rotationally symmetric manifolds; see \cite{Brendle2} for a precise statement.

We assume throughout that $g$ is a Riemannian metric on $\mathbb{R}^3$ such that
\begin{equation}
\label{expansiong}
g = (1 + |x|^{-1})^4 \sum_{i=1}^3 dx_i \otimes dx_i + T + o(|x|^{-2}) \text{ as } |x| \to \infty
\end{equation}
where the components of $T$ are homogeneous functions of degree $-2$. As in \cite{Huisken-Yau}, we require corresponding estimates for the first and second derivatives of the metric. The metric $g$ describes the asymptotically flat end of an initial data set whose mass is normalized to $2$. In fact, the first term on the right-hand side of (\ref{expansiong}) is the exact Schwarzschild metric. Our assumptions on the structure of the perturbation term $T$ are slightly stronger than in \cite{Huisken-Yau}. 

Let $\Sigma$ be a closed stable constant mean curvature surface in $(\mathbb{R}^3, g)$. Let $H_\Sigma$ denote the mean curvature of $\Sigma$ and let $\rho_\Sigma := \sup \{\rho > 0 : B_\rho(0) \cap \Sigma = \emptyset\}$. We say that $\Sigma$ is outlying if the compact region bounded by $\Sigma$ is disjoint from $B_{\rho_\Sigma}(0)$ and if the mean curvature of the coordinate spheres $\partial B_\rho(0)$ is positive for all $\rho \geq \rho_\Sigma$. By the maximum principle, $H_\Sigma > 0$ if $\Sigma$ is outlying.

Our first result shows that, in general, there may exist arbitrarily large outlying stable constant mean curvature spheres in an initial data set. This settles a question left open in the work of G.~Huisken and S.-T.~Yau; see \cite{Huisken-Yau}, p.~310.

\begin{theorem}
\label{thm1}
Let $a \in (0, \infty)$. There exists a metric $g$ as in (\ref{expansiong}) and a sequence of outlying stable constant mean curvature spheres $\Sigma^{(n)}$ such that $\rho_{\Sigma^{(n)}} \to \infty$, $H_{\Sigma^{(n)}} \to 0$, and such that $\lim_{n \to \infty} H_{\Sigma^{(n)}} \, \rho_{\Sigma^{(n)}} = 2 a \in (0, \infty)$.
\end{theorem}

On the other hand, we show that sequences of stable constant mean curvature spheres as in Theorem \ref{thm1} cannot exist when the scalar curvature of the metric $g$ is nonnegative.

\begin{theorem}
\label{thm2}
Let $a \in (0, \infty)$. Let $g$ be a metric as in (\ref{expansiong}) and let $R$ be its scalar curvature function. If there exists a sequence $\Sigma^{(n)}$ of outlying closed stable constant mean curvature surfaces such that $\rho_{\Sigma^{(n)}} \to \infty$ and $H_{\Sigma^{(n)}} \, \rho_{\Sigma^{(n)}} \to 2 a$ as $n \to \infty$, then $\liminf_{x \to \infty} |x|^4 \, R(x) < 0$.
\end{theorem}

By comparison, Theorem 1.6 in \cite{Eichmair-Metzger2} of J.~Metzger and the second-named author implies that if the scalar curvature is positive and if $\Sigma^{(n)}$ is a sequence of closed stable constant mean curvature surfaces whose areas diverge to infinity as $n \to \infty$, then $\rho_{\Sigma^{(n)}} \to \infty$.

Theorem \ref{thm2} may be rephrased as follows: Let $g$ be a metric as above such that $R \geq -o(|x|^{-4})$ where $R$ is the scalar curvature of $g$. Then, for every $\beta  \in (0, 1)$, we can find a real number $\rho$ with the property that there are no outlying stable constant mean curvature spheres $\Sigma$ with $\rho_\Sigma \geq \rho$ and $\beta \leq \rho_\Sigma H_{\Sigma} \leq \beta^{-1}$.

The uniqueness results in \cite{Huisken-Yau, Qing-Tian} for large stable constant mean curvature spheres that contain the center of the manifold  rely on a calculation of certain flux integrals. No assumption on the scalar curvature is necessary. By contrast, these flux integrals vanish in our situation. The proofs of our main results here are based on a more delicate analysis. Our starting point is a Lyapunov-Schmidt reduction. This method is inspired in part by the noncompactness results for the Yamabe equation in conformal geometry; see e.g. \cite{Ambrosetti} or \cite{Brendle1} for details. Using the implicit function theorem, we construct a family of surfaces $\Sigma_{(\xi,\lambda)}$ each of which is close to a large coordinate sphere centered at $\lambda \xi$ and of radius $\lambda$, whose enclosed volume is equal to $\frac{4\pi}{3} \lambda^3$, and whose mean curvature differs from a constant by first spherical harmonics. Moreover, we show that $\Sigma_{(\xi,\lambda)}$ has constant mean curvature with respect to $g$ if and only if $\xi$ is a critical point of the function $\xi \mapsto \mathscr{H}_g^2(\Sigma_{(\xi,\lambda)})$. On the other hand, we prove that
\[\mathscr{H}_g^2(\Sigma_{(\xi,\lambda)}) - 4\pi \, \lambda^2 = 2\pi \, F(\xi) + o(1),\]
where $F(\xi)$ is defined by
\begin{align}
\label{definition.of.F}
F(\xi)
&= -14 + 16 \, |\xi|^2 \log \frac{|\xi|^2-1}{|\xi|^2} + (15 \, |\xi| - |\xi|^{-1}) \, \log \frac{|\xi|+1}{|\xi|-1} \notag \\
&+ \frac{1}{4\pi} \, \bigg ( \int_{\partial B_1(\xi)} \text{\rm tr}_S(T) - 2 \int_{B_1(\xi)} \text{\rm tr}_{\mathbb{R}^3}(T) \bigg )
\end{align}
for $\xi \in \mathbb{R}^3 \setminus \bar{B}_1(0)$. It turns out that the last term is related to the scalar curvature of the background metric $g$ in a subtle way. If $R \geq -o(|x|^{-4})$, we are able to show that the radial derivative of the function $F$ is strictly positive. In particular, $F$ has no critical points in this case. On the other hand, if the hypothesis on the scalar curvature is dropped, we show that for a suitable choice of the perturbation term $T$ the function $F$ has a strict local minimum at some point $\xi \in \mathbb{R}^3$ with $|\xi|>1$.

\textbf{Acknowledgments.} The first-named author was supported in part by the U.S. National Science Foundation under grant DMS-1201924. He acknowledges the hospitality of the Department of Mathematics and Mathematical Statistics at Cambridge University, where part of this work was carried out. The second-named author was supported by the Swiss National Science Foundation under grant SNF 200021-140467. We are grateful to Professor Gerhard Huisken and Professor Jan Metzger for their interest and encouragement. Finally, we thank Otis Chodosh and the referee for their helpful remarks.

\section{Lyapunov-Schmidt reduction}

Let us fix a real number $\sigma > 1$, an integer $k \geq 1$, and a real number $\alpha \in (0, 1)$. Let $\mathscr{X}$ denote the space of Riemannian metrics on the ball $\bar{B}_\sigma(0) \subset \mathbb{R}^3$ with the $C^{k, \alpha}$ topology. Let $\mathscr{Y}$ denote the space of $C^{k+1, \alpha}$ functions on the unit sphere with center at the origin. The following consequence of the implicit function theorem is well known.

\begin{proposition}
\label{implicit.function.theorem}
If $g \in \mathscr{X}$ is sufficiently close to the Euclidean metric, then there exists a unique function $u \in \mathscr{Y}$ close to $0$ with the following properties:
\begin{itemize}
\item The spherical graph $\Sigma = \{(1+u(y)) \, y: y \in \partial B_1(0)\}$ encloses the volume $4 \pi / 3$ with respect to the metric $g$.
\item The function $u$ is orthogonal to the first spherical harmonics. The mean curvature of $\Sigma$ with respect to $g$ is equal to a constant plus a sum of first spherical harmonics.
\end{itemize}
We may write $u = \mathscr{G}(g)$ where $\mathscr{G}$ is a smooth map from a neighborhood of the Euclidean metric in $\mathscr{X}$ to a neighborhood of $0$ in $\mathscr{Y}$.
\end{proposition}

Let $g$ be a Riemannian metric of the form  (\ref{expansiong}). Fix a bounded open set $\Omega \subset \mathbb{R}^3$ with $\bar{\Omega} \cap \bar{B}_1(0) = \emptyset$. For $\xi \in \Omega$ and $\lambda>0$ large, we consider the coordinate sphere
\[S_{(\xi,\lambda)} = \{x \in \mathbb{R}^3:  |x-\lambda \, \xi| = \lambda\}.\]
The rescaled metric
\[\tilde g_{(\xi, \lambda)} = \lambda^{-2} \, \Phi_{(\xi,\lambda)}^* g\]
differs from the Euclidean metric by terms of order $O(\lambda^{-1})$ in $\mathscr{X}$. Here, $\Phi_{(\xi,\lambda)} : \mathbb{R}^3 \to \mathbb{R}^3$ is the map $y \mapsto \lambda (\xi+y)$. The following result follows from Proposition \ref{implicit.function.theorem} and a scaling argument.

\begin{proposition}
\label{deformation}
Suppose that $\lambda$ is sufficiently large. For every $\xi \in \Omega$ we can find a surface $\Sigma_{(\xi,\lambda)}$ with the following properties:
\begin{itemize}
\item $\Sigma_{(\xi,\varepsilon)}$ is a spherical  graph over the coordinate sphere $S_{(\xi,\lambda)}$. The height function $u_{(\xi,\lambda)}$ is orthogonal to the first spherical harmonics on $S_{(\xi,\lambda)}$ and we have the estimate
\[\sup_{S_{(\xi,\lambda)}} \left( |u_{(\xi,\lambda)}| + \lambda \, |D u_{(\xi,\lambda)}| + \lambda^2 \, |D^2 u_{(\xi,\lambda)}| \right) = O(1).\]
\item $\text{\rm vol}_g(\Sigma_{(\xi,\lambda)}) = 4 \pi \lambda^3 / 3$.
\item The mean curvature of $\Sigma_{(\xi,\lambda)}$ satisfies
\[H_{\Sigma_{(\xi,\lambda)}} - 2\lambda^{-1} = h_0(\xi,\lambda) + \sum_{i=1}^3 h_i(\xi,\lambda) \, y_i\]
where $h_i (\xi, \lambda) \in \mathbb{R}$ for $i = 0, 1, 2, 3$.
\end{itemize}
\end{proposition}

The enclosed volume of the coordinate sphere $S_{(\xi,\lambda)}$ with respect to the metric $g$ is given by
\begin{align*}
\text{\rm vol}_g(S_{(\xi,\lambda)})
&= \int_{\{x \in \mathbb{R}^3:  |x-\lambda \, \xi| \leq \lambda\}} (1 + |x|^{-1})^6 + O(\lambda) \\
&= \int_{\{x \in \mathbb{R}^3:  |x-\lambda \, \xi| \leq \lambda\}} (1 + 6 \, |x|^{-1}) + O(\lambda) \\
&= \frac{4\pi}{3} \, \lambda^3 \, (1 + 6 \, \lambda^{-1} \, |\xi|^{-1}) + O(\lambda).
\end{align*}
Using that $\text{\rm vol}_g(\Sigma_{(\xi,\lambda)}) = 4 \pi \lambda^3 / 3$ we find that
\begin{align*}
\int_{S_{(\xi,\lambda)}} u_{(\xi,\lambda)}
&= \text{\rm vol}_g(\Sigma_{(\xi,\lambda)}) - \text{\rm vol}_g(S_{(\xi,\lambda)}) + O(\lambda) \\
&= \frac{4\pi}{3} \, \lambda^3 - \text{\rm vol}_g(S_{(\xi,\lambda)}) + O(\lambda) \\
&= -8\pi \, \lambda^2 \, |\xi|^{-1} + O(\lambda).
\end{align*}
In other words, the mean value of $u_{(\xi,\lambda)}$ over the coordinate sphere $S_{(\xi,\lambda)}$ is equal to $-2 \, |\xi|^{-1} + O(\lambda^{-1})$.

The mean curvature $H_{S_{(\xi,\lambda)}}$ of the coordinate sphere $S_{(\xi,\lambda)}$ with respect to the metric $g$ satisfies
\[H_{S_{(\xi,\lambda)}} - 2\lambda^{-1} = O(\lambda^{-2}).\]
We will require a more precise estimate. Using the formula for the change of the mean curvature under a conformal change of the metric, we obtain that
\[H_{S_{(\xi,\lambda)}}-2\lambda^{-1} = -4\lambda^{-2} \, (|y+\xi|^{-1} + |y+\xi|^{-3} \, \langle y,y+\xi \rangle) + O(\lambda^{-3})\]
where $y = \lambda^{-1} \, x - \xi$. Since $|y|=1$ on $S_{(\xi,\lambda)}$ and $|\xi| > 1$, we may write
\[|y+\xi|^{-1} = \sum_{l \geq 0} |\xi|^{-l-1} \, P_l ( -\langle y,\xi \rangle / |\xi| ),\]
where the $P_l$'s are the Legendre polynomials. Differentiating this identity with respect to $\xi$ in radial direction gives
\[|y+\xi|^{-3} \, \langle \xi,y+\xi \rangle = \sum_{l \geq 0} (l+1) \, |\xi|^{-l-1} \, P_l ( -\langle y,\xi \rangle / |\xi| ).\]
Putting these facts together, we obtain
\begin{align*}
H_{S_{(\xi,\lambda)}}-2\lambda^{-1}
&= 4\lambda^{-2} \, (|y+\xi|^{-3} \, \langle \xi,y+\xi \rangle - 2 \, |y+\xi |^{-1}) + O(\lambda^{-3}) \\
&= 4\lambda^{-2} \sum_{l \geq 0} (l-1) \, |\xi|^{-l-1} \, P_l ( -\langle y,\xi \rangle / |\xi| ) + O(\lambda^{-3}).
\end{align*}
On the other hand, the mean curvature of the surface $\Sigma_{(\xi,\lambda)}$ with respect to the metric $g$ satisfies
\[H_{\Sigma_{(\xi,\lambda)}} = H_{S_{(\xi,\lambda)}} - \Delta_{S_{(\xi,\lambda)}} u_{(\xi,\lambda)} - 2\lambda^{-2} \, u_{(\xi,\lambda)} + O(\lambda^{-3}).\]
From this, we deduce that
\begin{align*}
\Delta_{S_{(\xi,\lambda)}} u_{(\xi,\lambda)} + 2\lambda^{-2} \, u_{(\xi,\lambda)}
&= H_{S_{(\xi,\lambda)}} - H_{\Sigma_{(\xi,\lambda)}} + O(\lambda^{-3}) \\
&= 4\lambda^{-2} \sum_{l \geq 0} (l-1) \, |\xi|^{-l-1} \, P_l ( -\langle y,\xi \rangle / |\xi| ) \\
&- h_0(\xi,\lambda) - \sum_{i=1}^3 h_i(\xi,\lambda) \, y_i + O(\lambda^{-3}).
\end{align*}
Using that $u_{(\xi,\lambda)}$ has mean $-2 \, |\xi|^{-1} + O(\lambda^{-1})$ and that it is orthogonal to the first spherical harmonics on $S_{(\xi,\lambda)}$, we obtain that
\[u_{(\xi,\lambda)} = -4 \sum_{l \geq 0, \, l \neq 1} \frac{1}{l+2} \, |\xi|^{-l-1} \, P_l ( -\langle y,\xi \rangle / |\xi| )+ O(\lambda^{-1})\]
and
\[|h_0(\xi,\lambda)| + \sum_{i=1}^3 |h_i(\xi,\lambda)| = O(\lambda^{-3}).\]
In particular, we have
\[H_{\Sigma_{(\xi,\lambda)}} - 2\lambda^{-1} = h_0(\xi,\lambda) + \sum_{i=1}^3 h_i(\xi,\lambda) \, y_i = O(\lambda^{-3}).\]
In the next step, we derive an asymptotic expansion for the surface area of $\Sigma_{(\xi,\varepsilon)}$.

\begin{proposition}
\label{expansion.for.area}
We have
\[\mathscr{H}_g^2(\Sigma_{(\xi,\lambda)}) - 4\pi \, \lambda^2  = 2\pi \, F(\xi) + o(1),\]
where $F(\xi)$ is defined by (\ref{definition.of.F}).
\end{proposition}

\textbf{Proof.} Let us define a functional $\mathscr{F}_\lambda$ on closed hypersurfaces $\Sigma$ by
\[\mathscr{F}_\lambda(\Sigma) = \mathscr{H}_g^2(\Sigma) - 2\lambda^{-1} \, \text{\rm vol}_g(\Sigma),\]
where the area and enclosed volume of $\Sigma$ are computed with respect to the metric $g$. Using the second variation formula for the functional $\mathscr{F}_\lambda$, we obtain
\begin{align*}
\mathscr{F}_\lambda(\Sigma_{(\xi,\lambda)})
&= \mathscr{F}_\lambda(S_{(\xi,\lambda)}) + \int_{S_{(\xi,\lambda)}} u_{(\xi,\lambda)} \, (H_{S_{(\xi,\lambda)}} - 2\lambda^{-1}) \\
&+ \frac{1}{2} \int_{S_{(\xi,\lambda)}} (|\nabla^{S_{(\xi,\lambda)}} u_{(\xi,\lambda)}|^2 - 2\lambda^{-2} \, u_{(\xi,\lambda)}^2) + O(\lambda^{-1}) \\
&= \mathscr{F}_\lambda(S_{(\xi,\lambda)}) - 16 \sum_{l \geq 0} \frac{l-1}{l+2} \, |\xi|^{-2l-2} \int_{\{|y|=1\}} P_l ( -\langle y,\xi \rangle / |\xi| )^2 \\
&+ 8 \sum_{l \geq 0} \frac{l-1}{l+2} \, |\xi|^{-2l-2} \int_{\{|y|=1\}} P_l ( -\langle y,\xi \rangle / |\xi| )^2 + O(\lambda^{-1}) \\
&= \mathscr{F}_\lambda(S_{(\xi,\lambda)}) - 16\pi \sum_{l \geq 0} \frac{l-1}{l+2} \, |\xi|^{-2l-2} \int_{-1}^1 P_l(-z)^2 \, dz + O(\lambda^{-1}) \\
&= \mathscr{F}_\lambda(S_{(\xi,\lambda)}) - 32\pi \sum_{l \geq 0} \frac{l-1}{(l+2)(2l+1)} \, |\xi|^{-2l-2} + O(\lambda^{-1}).
\end{align*}
Note that
\[\sum_{l \geq 0} \frac{1}{l+2} \, |\xi|^{-2l-2} = -1 - |\xi|^2 \log \frac{|\xi|^2-1}{|\xi|^2}\]
and
\[\sum_{l \geq 0} \frac{1}{2l+1} \, |\xi|^{-2l-2} = \frac{1}{2 \, |\xi|} \log \frac{|\xi|+1}{|\xi|-1}.\]
Subtracting the second identity from the first gives
\[\sum_{l \geq 0} \frac{l-1}{(l+2)(2l+1)} \, |\xi|^{-2l-2} = -1 - |\xi|^2 \log \frac{|\xi|^2-1}{|\xi|^2} - \frac{1}{2 \, |\xi|} \log \frac{|\xi|+1}{|\xi|-1}.\]
Thus,
\[\mathscr{F}_\lambda(\Sigma_{(\xi,\lambda)}) = \mathscr{F}_\lambda(S_{(\xi,\lambda)}) + 32\pi \, \Big ( 1 + |\xi|^2 \log \frac{|\xi|^2-1}{|\xi|^2} + \frac{1}{2 \, |\xi|} \log \frac{|\xi|+1}{|\xi|-1} \Big ) + O(\lambda^{-1}).\]
We next analyze the term $\mathscr{F}_\lambda(S_{(\xi,\lambda)})$. Using the identity
\[\int_{S_{(\xi,\lambda)}} |x|^{-2} = 2\pi \, |\xi|^{-1} \, \log \frac{|\xi|+1}{|\xi|-1}\]
we obtain
\begin{align*}
\mathscr{H}_g^2(S_{(\xi,\lambda)})
&= \int_{S_{(\xi,\lambda)}} (1 + |x|^{-1})^4 + \frac{1}{2} \int_{S_{(\xi,\lambda)}} \text{\rm tr}_{S_{(\xi,\lambda)}}(T) + o(1) \\
&= \int_{S_{(\xi,\lambda)}} (1 + 4 \, |x|^{-1} + 6 \, |x|^{-2}) + \frac{1}{2} \int_{S_{(\xi,\lambda)}} \text{\rm tr}_{S_{(\xi,\lambda)}}(T) + o(1) \\
&= 4\pi \, \lambda^2 + 16\pi \, \lambda \, |\xi|^{-1} \\
&+ 12\pi \, |\xi|^{-1} \, \log \frac{|\xi|+1}{|\xi|-1} + \frac{1}{2} \int_{\partial B_1(\xi)} \text{\rm tr}_S(T) + o(1).
\end{align*}
Similarly, the identity
\[\int_{B_{(\xi,\lambda)}} |x|^{-2} = 2\pi \, \lambda \, \Big ( 1 - \frac{|\xi|^2-1}{2 \, |\xi|} \, \log \frac{|\xi|+1}{|\xi|-1} \Big )\]
gives
\begin{align*}
\text{\rm vol}_g(S_{(\xi,\lambda)})
&= \int_{B_{(\xi,\lambda)}} (1 + |x|^{-1})^6 + \frac{1}{2} \int_{B_{(\xi,\lambda)}} \text{\rm tr}_{\mathbb{R}^3}(T) + o(\lambda) \\
&= \int_{B_{(\xi,\lambda)}} (1 + 6 \, |x|^{-1} + 15 \, |x|^{-2}) + \frac{1}{2} \int_{B_{(\xi,\lambda)}} \text{\rm tr}_{\mathbb{R}^3}(T) + o(\lambda) \\
&= \frac{4\pi}{3} \, \lambda^3 + 8\pi \, \lambda^2 \, |\xi|^{-1} + 30\pi \, \lambda \, \Big ( 1 - \frac{|\xi|^2-1}{2 \, |\xi|} \, \log \frac{|\xi|+1}{|\xi|-1} \Big ) \\
&+ \frac{1}{2} \, \lambda \int_{B_1(\xi)} \text{\rm tr}_{\mathbb{R}^3}(T) + o(\lambda).
\end{align*}
Putting these facts together, we conclude that
\begin{align*}
\mathscr{F}_\lambda(S_{(\xi,\lambda)})
&= \frac{4\pi}{3} \, \lambda^2 + 12\pi \, |\xi|^{-1} \, \log \frac{|\xi|+1}{|\xi|-1} - 60\pi \, \Big ( 1 - \frac{|\xi|^2-1}{2 \, |\xi|} \, \log \frac{|\xi|+1}{|\xi|-1} \Big ) \\
&+ \frac{1}{2} \int_{\partial B_1(\xi)} \text{\rm tr}_S(T) - \int_{B_1(\xi)} \text{\rm tr}_{\mathbb{R}^3}(T) + o(1).
\end{align*}
This implies
\begin{align*}
\mathscr{F}_\lambda(\Sigma_{(\xi,\lambda)})
&= \frac{4\pi}{3} \, \lambda^2 + 12\pi \, |\xi|^{-1} \, \log \frac{|\xi|+1}{|\xi|-1} - 60\pi \, \Big ( 1 - \frac{|\xi|^2-1}{2 \, |\xi|} \, \log \frac{|\xi|+1}{|\xi|-1} \Big ) \\
&+ 32\pi \Big ( 1 + |\xi|^2 \log \frac{|\xi|^2-1}{|\xi|^2} + \frac{1}{2 \, |\xi|} \log \frac{|\xi|+1}{|\xi|-1} \Big ) \\
&+ \frac{1}{2} \int_{\partial B_1(\xi)} \text{\rm tr}_S(T) - \int_{B_1(\xi)} \text{\rm tr}_{\mathbb{R}^3}(T) + o(1) \\
&= \frac{4\pi}{3} \, \lambda^2 + 2\pi \, F(\xi) + o(1).
\end{align*}
Since $\mathscr{F}_\lambda(\Sigma_{(\xi,\lambda)}) = \mathscr{H}_g^2(\Sigma_{(\xi,\lambda)}) - \frac{8\pi}{3} \, \lambda^2$, we obtain
\[\mathscr{H}_g^2(\Sigma_{(\xi,\lambda)}) = 4\pi \, \lambda^2 + 2\pi \, F(\xi) + o(1).\]
This completes the proof of Proposition \ref{expansion.for.area}. \\

Finally, we need the following technical estimate:

\begin{proposition}
\label{technical}
Let $\Omega \subset \mathbb{R}^3$ be a bounded open set with $\bar{\Omega} \cap \bar{B}_1(0) = \emptyset$. Then
\[\sup_\Omega \max_{1 \leq i \leq 3} \Big | \frac{\partial}{\partial \xi_i} \mathscr{H}_g^2(\Sigma_{(\xi,\lambda)}) \Big | \leq C  \text{ and }    \sup_\Omega \max_{1 \leq i,j \leq 3} \Big | \frac{\partial^2}{\partial \xi_i \, \partial \xi_j} \mathscr{H}_g^2(\Sigma_{(\xi,\lambda)}) \Big | \leq C,\]
where the constant $C$ is independent of $\lambda$.
\end{proposition}

\textbf{Proof.}
The rescaled surface $\tilde{\Sigma}_{(\xi,\lambda)} = \lambda^{-1} \, \Sigma_{(\xi,\lambda)} - \xi$ is a spherical graph over the unit sphere with center at the origin. Its height function $\tilde{u}_{(\xi,\lambda)}$ can be expressed as $\tilde{u}_{(\xi,\lambda)} = \mathscr{G}(\tilde{g}_{(\xi,\lambda)})$ where $\mathscr{G}$ is a $C^2$ map from a neighborhood of the Euclidean metric in $\mathscr{X}$ to a neighborhood of $0$ in $\mathscr{Y}$. Since
\[\Big \| \frac{\partial}{\partial \xi_i} \tilde{g}_{(\xi,\lambda)} \Big \|_{\mathscr{X}} \leq O(\lambda^{-1}) \quad \text{ and } \quad   \Big \| \frac{\partial^2}{\partial \xi_i \, \partial \xi_j} \tilde{g}_{(\xi,\lambda)} \Big \|_{\mathscr{X}} \leq O(\lambda^{-1}),\]
we conclude that
\[\Big \| \frac{\partial}{\partial \xi_i} \tilde{u}_{(\xi,\lambda)} \Big \|_{\mathscr{Y}} \leq O(\lambda^{-1})\quad \text{ and } \quad   \Big \| \frac{\partial^2}{\partial \xi_i \, \partial \xi_j} \tilde{u}_{(\xi,\lambda)} \Big \|_{\mathscr{Y}} \leq O(\lambda^{-1}).\]
Fix a point $\bar{\xi} \in \Omega$. Using the formulae for the first variation of area and volume, we obtain
\begin{align*}
&\frac{\partial}{\partial \xi_i} \mathscr{H}_{\tilde{g}_{(\xi,\lambda)}}^2(\tilde{\Sigma}_{(\xi,\lambda)}) \Big |_{\xi=\bar{\xi}} \\
&= \frac{\partial}{\partial \xi_i} \Big ( \mathscr{H}_{\tilde{g}_{(\xi,\lambda)}}^2(\tilde{\Sigma}_{(\xi,\lambda)}) - 2 \, \text{\rm vol}_{\tilde{g}_{(\xi,\lambda)}}(\tilde{\Sigma}_{(\xi,\lambda)}) \Big ) \Big |_{\xi=\bar{\xi}} \\
&= \frac{\partial}{\partial \xi_i} \Big ( \mathscr{H}_{\tilde{g}_{(\xi,\lambda)}}^2(\tilde{\Sigma}_{(\bar{\xi},\lambda)}) - 2 \, \text{\rm vol}_{\tilde{g}_{({\xi},\lambda)}}(\tilde{\Sigma}_{(\bar{\xi},\lambda)}) \Big ) \Big |_{\xi=\bar{\xi}} + O(\lambda^{-2}) \\
&= - 4\lambda^{-1} \int_{\partial B_1(0)} |\xi+y|^{-3} \, (\xi_i+y_i) \\
& +12\lambda^{-1} \int_{B_1(0)} |\xi+y|^{-3} \, (\xi_i+y_i) + O(\lambda^{-2}) \\
&= O(\lambda^{-2}).
\end{align*}
Similarly, using the formula for the second variation of area and volume, we obtain
\begin{align*}
&\frac{\partial^2}{\partial \xi_i \, \partial \xi_j} \mathscr{H}_{\tilde{g}_{(\xi,\lambda)}}^2(\tilde{\Sigma}_{(\xi,\lambda)}) \Big |_{\xi=\bar{\xi}} \\
&= \frac{\partial^2}{\partial \xi_i \, \partial \xi_j} \Big ( \mathscr{H}_{\tilde{g}_{(\xi,\lambda)}}^2(\tilde{\Sigma}_{(\xi,\lambda)}) - 2 \, \text{\rm vol}_{\tilde{g}_{(\xi,\lambda)}}(\tilde{\Sigma}_{(\xi,\lambda)}) \Big ) \Big |_{\xi=\bar{\xi}} \\
&= \frac{\partial^2}{\partial \xi_i \, \partial \xi_j} \Big ( \mathscr{H}_{\tilde{g}_{(\xi,\lambda)}}^2(\tilde{\Sigma}_{(\bar{\xi},\lambda)}) - 2 \, \text{\rm vol}_{\tilde{g}_{({\xi},\lambda)}}(\tilde{\Sigma}_{(\bar{\xi},\lambda)}) \Big ) \Big |_{\xi=\bar{\xi}} + O(\lambda^{-2}) \\
&= 4\lambda^{-1} \int_{\partial B_1(0)} (|\xi+y|^{-3} \, \delta_{ij} - 3 \, |\xi+y|^{-5} \, (\xi_i+y_i) \, (\xi_j+y_j)) \\
&- 12\lambda^{-1} \int_{B_1(0)} (|\xi+y|^{-3} \, \delta_{ij} - 3 \, |\xi+y|^{-5} \, (\xi_i+y_i) \, (\xi_j+y_j)) + O(\lambda^{-2}) \\
&= O(\lambda^{-2}).
\end{align*}

\section{Perturbations with $R \geq -o(|x|^{-4})$}

Throughout this section, we will assume that $R \geq -o(|x|^{-4})$. It follows from (\ref{expansiong}) that the scalar curvature of $g$ is given by
\[R = \sum_{i,j=1}^3 (D_i D_j T_{ij} - D_i D_i T_{jj}) + O(|x|^{-5}).\]
The condition $R \geq -o(|x|^{-4})$ is equivalent to the inequality
\begin{equation}
\label{T}
\sum_{i,j=1}^3 (D_i D_j T_{ij} - D_i D_i T_{jj}) \geq 0.
\end{equation}

\begin{proposition}
\label{nonnegative.scalar.curvature}
Let $\xi \in \mathbb{R}^3$ be such that $|\xi| > 1$. If $T$ satisfies (\ref{T}), then
\[\frac{d}{ds} \bigg ( \int_{\partial B_1(s\xi)} \text{\rm tr}_S(T) - 2 \int_{B_1(s\xi)} \text{\rm tr}_{\mathbb{R}^3}(T) \bigg ) \bigg |_{s=1} \geq 0.\]
\end{proposition}

\textbf{Proof.}
Let
\[G(s) = \int_{\partial B_1(s\xi)} \text{\rm tr}_S(T) - 2 \int_{B_1(s\xi)} \text{\rm tr}_{\mathbb{R}^3}(T)\]
and
\[K(s) = \int_{B_1(s\xi)} \sum_{i,j=1}^3 (D_i D_j T_{ij} - D_i D_i T_{jj}).\]
By assumption, $K(s)$ is nonnegative, and we have
\[K(s) \leq \int_{B_s(s\xi)} \sum_{i,j=1}^3 (D_i D_j T_{ij} - D_i D_i T_{jj}) = K(1) / s \]
for all $s \geq 1$.

Using the divergence theorem, we obtain
\begin{align*}
K(s) &= \int_{\partial B_1(s\xi)} \sum_{a=1}^2 \big ( (D_{e_a} T)(e_a,\nu) - (D_\nu T)(e_a,e_a) \big ) \\
&= \int_{\partial B_1(s\xi)} \bigg ( 2 \, T(\nu,\nu) - \sum_{a=1}^2 T(e_a,e_a) - \sum_{a=1}^2 (D_\nu T)(e_a,e_a) \bigg ) \\
&= \int_{\partial B_1(s\xi)} \big ( 2 \, \text{\rm tr}_{\mathbb{R}^3}(T) - 3 \, \text{\rm tr}_S(T) - \text{\rm tr}_S(D_\nu T) \big ).
\end{align*}
We next observe that
\[\frac{d}{ds} \bigg ( \int_{\partial B_1(s\xi)} \text{\rm tr}_S(T) \bigg ) = \int_{\partial B_1(s\xi)} \text{\rm tr}_S(D_\xi T)\]
and
\[\frac{d}{ds} \bigg ( \int_{B_1(s\xi)} \text{\rm tr}_{\mathbb{R}^3}(T) \bigg ) = \int_{B_1(s\xi)} \text{\rm tr}_{\mathbb{R}^3}(D_\xi T) = \int_{\partial B_1(s\xi)} \text{\rm tr}_{\mathbb{R}^3}(T) \, \langle \xi,\nu \rangle.\]
Since $T$ is homogeneous of degree $-2$, the radial derivative of $T$ satisfies $D_x T = -2 \, T$ by Euler's theorem. Putting these facts together, we obtain
\begin{align*}
s \, G'(s)
&= \int_{\partial B_1(s\xi)} \text{\rm tr}_S(D_{s\xi} T) - 2 \int_{\partial B_1(s\xi)} \text{\rm tr}_{\mathbb{R}^3}(T) \, \langle s\xi,\nu \rangle \\
&= \int_{\partial B_1(s\xi)} \text{\rm tr}_S(D_x T) - 2 \int_{\partial B_1(s\xi)} \text{\rm tr}_{\mathbb{R}^3}(T) \, \langle x,\nu \rangle \\
&- \int_{\partial B_1(s\xi)} \text{\rm tr}_S(D_\nu T) + 2 \int_{\partial B_1(s\xi)} \text{\rm tr}_{\mathbb{R}^3}(T) \\
&= \int_{\partial B_1(s\xi)} \text{\rm tr}_S(D_x T) - 2 \int_{B_1(s\xi)} (\text{\rm tr}_{\mathbb{R}^3}(D_x T) + 3 \, \text{\rm tr}_{\mathbb{R}^3}(T)) \\
&- \int_{\partial B_1(s\xi)} \text{\rm tr}_S(D_\nu T) + 2 \int_{\partial B_1(s\xi)} \text{\rm tr}_{\mathbb{R}^3}(T) \\
&= -2 \int_{\partial B_1(s\xi)} \text{\rm tr}_S(T) - 2 \int_{B_1(s\xi)} \text{\rm tr}_{\mathbb{R}^3}(T) \\
&- \int_{\partial B_1(s\xi)} \text{\rm tr}_S(D_\nu T) + 2 \int_{\partial B_1(s\xi)} \text{\rm tr}_{\mathbb{R}^3}(T) \\
&= G(s) + K(s)
\end{align*}
for all $s \geq 1$. Here, $\nu = x - s\xi$ denotes the outward-pointing unit normal vector along the sphere $\partial B_1(s\xi)$.

Since $K(s) \leq K(1)/s$, it follows that the function $G(s)/2 + K(1) / ( 2s^2)$ is monotone decreasing in $s$. Since $G(s) /s  +  K(1) / (2 s^2) \to 0$ as $s \to \infty$, we conclude that
\[G(1) + K(1)/2 \geq 0.\]
Since $K(1) \geq 0$, it follows that
\[G'(1) = G(1) + K(1) \geq 0,\]
as claimed. \\

\begin{corollary}
\label{radial.derivative}
Assume that $T$ satisfies (\ref{T}). Then we have
\[\frac{d}{ds} F(s\xi) \Big |_{s=1} \geq 32 \, |\xi|^2 \log \frac{|\xi|^2-1}{|\xi|^2} + (15 \, |\xi|+|\xi|^{-1}) \log \frac{|\xi|+1}{|\xi|-1} + 2 \, \frac{|\xi|^2+1}{|\xi|^2-1},\]
where $F(\xi)$ is defined by (\ref{definition.of.F}).
\end{corollary}

\textbf{Proof.}
From the definition of $F$, we have that
\begin{align*}
\frac{d}{ds} F(s\xi) \Big |_{s=1}
&= 32 \, |\xi|^2 \log \frac{|\xi|^2-1}{|\xi|^2} + (15 \, |\xi|+|\xi|^{-1}) \log \frac{|\xi|+1}{|\xi|-1} + 2 \, \frac{|\xi|^2+1}{|\xi|^2-1} \\
&+ \frac{1}{4\pi} \, \frac{d}{ds} \bigg ( \int_{\partial B_1(s\xi)} \text{\rm tr}_S(T) - 2 \int_{B_1(s\xi)} \text{\rm tr}_{\mathbb{R}^3}(T) \bigg ) \bigg |_{s=1}.
\end{align*}
The assertion now follows from Proposition \ref{nonnegative.scalar.curvature}.

\section{Proof of Theorem \ref{thm2}}

We argue by contradiction. Suppose that $\lim_{x \to \infty} |x|^4 \, R(x) \geq 0$. Let $\Sigma^{(n)}$ be a sequence of outlying stable constant mean curvature surfaces such that $\rho_{\Sigma^{(n)}} \to \infty$ and $H_{\Sigma^{(n)}} \to 0$. Let $r_n:= 2 /  H_{\Sigma^{(n)}}$. By assumption, we have that
\[\lim_{n \to \infty}  \rho_{\Sigma^{(n)}} / r_n = \lim_{n \to \infty} H_{\Sigma^{(n)}} \, \rho_{\Sigma^{(n)}} /2= a \in (0,\infty).\]
Consider the metrics $g^{(n)} := r_n^{-2} \, \varphi_n^* g$ where $\varphi_n : \mathbb{R}^3 \to \mathbb{R}^3$ is the homothety $x \mapsto r_n x$. The rescaled metrics $g^{(n)}$ converge to the Euclidean metric $\sum_{i=1}^3 dx_i \otimes dx_i$ away from the origin. The mean curvature of the surfaces $r_n^{-1} \, \Sigma^{(n)}$ with respect to $g^{(n)}$ is equal to $2$. The second fundamental forms of these surfaces and their derivatives are uniformly bounded by Proposition 2.2 in \cite{Eichmair-Metzger2}. Using compactness arguments for immersions with bounded geometry as in \cite{Cooper} we may extract a subsequence that converges to a stable constant mean curvature immersion with mean curvature $2$ in $\mathbb{R}^3$. Such an immersion is a unit sphere by results of B.~Palmer \cite{Palmer}; see also \cite{daSilveira} and \cite{Lopez-Ros}. To summarize, a subsequence of the rescaled surfaces $r_n^{-1} \, \Sigma^{(n)}$ converges to a unit sphere $\partial B_1(\bar{\xi})$ where $|\bar{\xi}| = 1 + a$. It follows that $\Sigma^{(n)}$ is a perturbation of a coordinate sphere when $n$ is sufficiently large and thus arises in our Lyapunov-Schmidt reduction. More precisely, we may write $\Sigma^{(n)} = \Sigma_{(\xi_n,\lambda_n)}$ where $\lim_{n \to \infty} \xi_n \to \bar{\xi}$ and $\lim_{n \to \infty} r_n / \lambda_n \to 1$.

We claim that this setup leads to a contradiction. To see this, let $\Omega$ denote the open ball centered at $\bar{\xi}$ of radius $a/2$. Clearly, $\xi_n$ lies in the interior of $\Omega$ if $n$ is sufficiently large, and we have $\bar{\Omega} \cap \bar{B}_1(0) = \emptyset$. We define a function $F_\lambda: \Omega \to \mathbb{R}$ by
\[F_\lambda(\xi) = \frac{1}{2 \pi} \left( \mathscr{H}_g^2(\Sigma_{(\xi,\lambda)}) - 4 \pi \lambda^2 \right).\]
Since $\Sigma^{(n)}$ has constant mean curvature, it follows that $\xi_n$ is a critical point of the function $F_{\lambda_n}$.

By Proposition \ref{expansion.for.area}, we have $\|F_\lambda - F\|_{C^0(\Omega)} \to 0$ as $\lambda \to \infty$. Moreover, Proposition \ref{technical} implies that $\|F_\lambda\|_{C^2(\Omega)} \leq C$ for some constant $C$ that is independent of $\lambda > 0$. Putting these facts together, we conclude that $\|F_\lambda - F\|_{C^1(\Omega)} \to 0$ as $\lambda \to \infty$. Using Corollary \ref{radial.derivative}, we obtain that
\[\frac{d}{ds} F_\lambda(s\xi) \Big |_{s=1} \geq 32 \, |\xi|^2 \log \frac{|\xi|^2-1}{|\xi|^2} + (15 \, |\xi|+|\xi|^{-1}) \log \frac{|\xi|+1}{|\xi|-1} + 2 \, \frac{|\xi|^2+1}{|\xi|^2-1} - o(1)\]
for all $\xi \in \Omega$. On the other hand, it is elementary to see that
\[32 \, |\xi|^2 \log \frac{|\xi|^2-1}{|\xi|^2} + (15 \, |\xi|+|\xi|^{-1}) \log \frac{|\xi|+1}{|\xi|-1} + 2 \, \frac{|\xi|^2+1}{|\xi|^2-1} > 0\]
for all $\xi \in \mathbb{R}^3$ with $|\xi| > 1$.\footnote{Note that if the left-hand side were negative for some $\xi >1$, then we could turn the argument around and conclude the existence of closed constant mean curvature surfaces in the Schwarzschild manifold that are disjoint from the horizon. This would contradict results in \cite{Brendle2}.} Thus
\[\inf_{\xi \in \Omega} \Big ( 32 \, |\xi|^2 \log \frac{|\xi|^2-1}{|\xi|^2} + (15 \, |\xi|+|\xi|^{-1}) \log \frac{|\xi|+1}{|\xi|-1} + 2 \, \frac{|\xi|^2+1}{|\xi|^2-1} \Big ) > 0.\]
Hence, if we choose $\lambda>0$ sufficiently large, we have that
\[\frac{d}{ds} F_\lambda(s\xi) \Big |_{s=1} > 0\]
for all points $\xi \in \Omega$. In particular, if $\lambda$ is large enough, then the function $F_\lambda$ has no critical points in $\Omega$. This contradicts the fact that $\xi_n$ is a critical point of the function $F_{\lambda_n}$.

\section{Proof of Theorem \ref{thm1}}

We now construct a metric $g$ that admits a one-parameter family of outlying stable constant mean curvature spheres. Throughout this section, we let $v = (0,0,1) \in \mathbb{R}^3$ and we consider the function
\[I(s) = \int_{\partial B_1(sv)} |x|^{-2} \, (1 - 3 \, (x_3-s))^2\]
for $s>1$. Let us fix a real number $s_0 \in [2,\infty)$ such that $I'(s_0) \neq 0$.

\begin{lemma}
\label{a}
There exists a smooth function $\psi: \mathbb{R} \to \mathbb{R}$ such that the function
\[ \xi \mapsto  \int_{\partial B_1(\xi)} |x|^{-2} \, \psi(x_3/|x|) \, (1 - 3 \, (x_3-\xi_3)^2),\]
has a strict local minimum at some point $\xi \in \mathbb{R}^3$ with $|\xi|>1$.
\end{lemma}

\textbf{Proof.}
Let $t_0 = \sqrt{1- s_0^{-2} }$. Consider the smooth function
\[\varphi_k(t) = \begin{cases} e^{-\frac{4}{k \, (t_0-t)+1}} & \text{\rm if $t < t_0+\frac{1}{k}$} \\ 0 & \text{\rm if $t \geq t_0+\frac{1}{k}$}. \end{cases}\]
Note that $\sup_{t \geq t_0} |\varphi_k(t)| \leq 1$ and $\sup_{t \geq t_0} |\varphi_k'(t)| = O(k)$. We also define the function
\[J_k(s) = \int_{\partial B_1(sv)} |x|^{-2} \, \varphi_k(x_3 / |x|) \, (1 - 3 \, (x_3-s)^2).\]
Note that $x_3 / |x| \geq t_0$ for all $x \in B_1(s_0 v)$. It follows that $|J_k(s_0)| = O(1)$ and $|J_k'(s_0)| = O(k)$. Hence, if we define
\[a_k := \frac{J_k'(s_0)}{I'(s_0)},\]
then $|a_k| = O(k)$.

We next define a smooth function $\psi_k$ by $\psi_k(t) = \varphi_k(t) - a_k$. Moreover, we define
\[Q_k(\xi) = \int_{\partial B_1(\xi)} |x|^{-2} \, \psi_k(x_3 / |x|) \, (1 - 3 \, (x_3-\xi_3)^2).\]
Clearly, $Q_k(sv) = J_k(s) - a_k \, I(s)$ for all $s$. This implies that
\[\frac{d}{ds} Q_k(sv) \big |_{s=s_0} = J_k'(s_0) - a_k \, I'(s_0) = 0\]
by definition of $a_k$. It follows that the point $s_0  v$ is a critical point of the function $Q_k$.

We claim that the function $Q_k$ has a strict local minimum at the point $s_0 v$ if $k$ is sufficiently large. To see this, we examine the Hessian of the function $Q_k$ at the point $s_0 v$. For abbreviation, let
\[U_k := \{x \in \partial B_1(s_0 v): x_3 / |x| \in [t_0,t_0+1/k]\}\]
and
\[V_k := \{x \in \partial B_1(s_0 v): x_3 / |x| \in [t_0, t_0+1/(2k)]\}.\]
We have that
\begin{align*}
&\frac{\partial^2}{\partial \xi_i \, \partial \xi_j} Q_k \Big |_{\xi=s_0 v} \\
&= \int_{U_k} |x|^{-8} \, \varphi_k''(x_3 / |x|) \, (|x|^2 \, \delta_{i3} - x_i \, x_3) \, (|x|^2 \, \delta_{j3} - x_j \, x_3) \, (1 - 3 \, (x_3-s_0)^2) + O(k).
\end{align*}
For every point $x \in U_k$ we have $1 - 3 \, (x_3-s_0)^2 \geq 1/4$ and $\varphi_k''(x_3 / |x|) \geq 0$. Moreover, we have $\varphi_k''( x_3 / |x|) \geq c \, k^2$ for all $x \in V_k$. The measure of $V_k$ is bounded below by $c \, k^{-1/2 }$. It follows that
\begin{align*}
\frac{\partial^2}{\partial \xi_3^2} Q_k \Big |_{\xi=s_0 v}
&= \int_{U_k} |x|^{-8} \, \varphi_k''(x_3 / |x|) \, (x_1^2+x_2^2)^2 \, (1 - 3 \, (x_3-s_0)^2) + O(k) \\
&\geq \int_{V_k} |x|^{-8} \, \varphi_k''(x_3 / |x|) \, (x_1^2+x_2^2)^2 \, (1 - 3 \, (x_3-s_0)^2) + O(k) \\
&\geq c \, k^{3/2} - O(k)
\end{align*}
for some $c >0$. A similar calculation gives that
\begin{align*}
&\frac{\partial^2}{\partial \xi_1^2} Q_k \Big |_{\xi=s_0 v} = \frac{\partial^2}{\partial \xi_2^2} Q_k \Big |_{\xi=s_0 v} \\
&= \int_{U_k} |x|^{-8} \, \varphi_k''(x_3 / |x|) \, (x_1^2+x_2^2) \, x_3^2 \, (1 - 3 \, (x_3-s_0)^2) /2 + O(k) \\
&\geq \int_{V_k} |x|^{-8} \, \varphi_k''(x_3 / |x|) \, (x_1^2+x_2^2) \, x_3^2 \, (1 - 3 \, (x_3-s_0)^2) /2 + O(k) \\
&\geq c \, k^{3/2} - O(k)
\end{align*}
for some $c>0$. Moreover, since $Q_k$ is axially symmetric, we know that the Hessian of $Q_k$ at the point $s_0 v$ is a diagonal matrix. It follows that the Hessian of $Q_k$ at $s_0 v$ is positive definite if $k$ is sufficiently large.  \\

\begin{lemma}
\label{b}
There exists a smooth function $\psi: \mathbb{R} \to \mathbb{R}$ such that the function
\begin{align*}
\xi \mapsto &-28\pi + 32\pi \, |\xi|^2 \log \frac{|\xi|^2-1}{|\xi|^2} + 2\pi \, (15 \, |\xi| - |\xi|^{-1}) \, \log \frac{|\xi|+1}{|\xi|-1} \\
&+ \int_{\partial B_1(\xi)} |x|^{-2} \, \psi(x_3 / |x|) \, (1 - 3 \, (x_3-\xi_3)^2)
\end{align*}
has a strict local minimum at some point $\xi$ with $|\xi|>1$.
\end{lemma}

\textbf{Proof.}
Multiply the function $\psi$ from Lemma \ref{a} by a large positive constant.  \\

Let $\psi$ be as in Lemma \ref{b}. We define a Riemannian metric $g$ by
\[g = (1 + |x|^{-1})^4 \sum_{i=1}^3 dx_i \otimes dx_i + T\]
where
\[T = -2 \, |x|^{-2} \, \psi(x_3 / |x|) \, (dx_1 \otimes dx_1 + dx_2 \otimes dx_2 - 2 \, dx_3 \otimes dx_3).\]
Note that the components of $T$ are homogeneous functions of degree $-2$ and that $\text{\rm tr}_{\mathbb{R}^3}(T) = 0$. Using Proposition \ref{expansion.for.area}, we find that
\begin{align*}
&\mathscr{H}_g^2(\Sigma_{(\xi,\lambda)}) - 4\pi \, \lambda^2 \\
&= -28\pi + 32\pi \, |\xi|^2 \log \frac{|\xi|^2-1}{|\xi|^2} + 2\pi \, (15 \, |\xi| - |\xi|^{-1}) \, \log \frac{|\xi|+1}{|\xi|-1} \\
&+ \frac{1}{2} \int_{\partial B_1(\xi)} \text{\rm tr}_S(T) + o(1) \\
&= -28\pi + 32\pi \, |\xi|^2 \log \frac{|\xi|^2-1}{|\xi|^2} + 2\pi \, (15 \, |\xi| - |\xi|^{-1}) \, \log \frac{|\xi|+1}{|\xi|-1} \\
&+ \int_{\partial B_1(\xi)} |x|^{-2} \, \psi(x_3 / |x|) \, (1 - 3 \, (x_3-\xi_3)^2) + o(1).
\end{align*}
Hence, if $\lambda$ is sufficiently large, the function
\[\xi \mapsto \mathscr{H}_g^2(\Sigma_{(\xi,\lambda)}) - 4\pi \, \lambda^2\]
attains a local minimum at some point $\xi$ with $|\xi|>1$. Consequently, the surface $\Sigma_{(\xi, \lambda)}$ corresponding to this particular value of $\xi$ is an outlying stable constant mean curvature sphere. This finishes the proof of Theorem \ref{thm1}.


\begin{thebibliography}{99}
\bibitem{Alexandrov}
A.D.~Alexandrov, \textit{Uniqueness theorems for surfaces in the large I,} Vesnik Leningrad Univ. 11, 5--17 (1956)

\bibitem{Ambrosetti}
A.~Ambrosetti, \textit{Multiplicity results for the Yamabe problem on $S^n$,} Proc. Natl. Acad. Sci. USA 99, 15252--15256 (2002)

\bibitem{Bray-thesis}
H.~Bray, \textit{The Penrose inequality in general relativity and volume comparison theorems involving scalar curvature,} PhD thesis, Stanford University (1997)

\bibitem{Brendle1}
S.~Brendle, \textit{Blow-up phenomena for the Yamabe equation,} J. Amer. Math. Soc. 21, 951--979 (2008)

\bibitem{Brendle2}
S.~Brendle, \textit{Constant mean curvature surfaces in warped product manifolds,} Publ. Math. IH\'ES 117, 247--269 (2013)

\bibitem{Brendle-Eichmair}
S.~Brendle and M.~Eichmair, \textit{Isoperimetric and Weingarten surfaces in the Schwarzschild manifold,} J. Diff. Geom. 94, no. 3, 387--407 (2013)

\bibitem{Christodoulou-Yau}
D.~Christodoulou and S.-T.~Yau, \textit{Some remarks on the quasi-local mass,} Mathematics and general relativity (Santa Cruz, 1986), Contemporary Mathematics volume 71, pp.~9--14, Amer. Math. Soc., Providence RI (1986)

\bibitem{Cooper}
A.~Cooper, \textit{A Compactness Theorem for the Second Fundamental Form,} arXiv:1006.5697 [math.DG]

\bibitem{daSilveira} A.~Da Silveira, \textit{Stability of complete noncompact surfaces with constant mean curvature,} Math. Ann. 277 (1987)

\bibitem{Eichmair-Metzger1}
M.~Eichmair and J.~Metzger, \textit{Large isoperimetric surfaces in initial data sets,} J. Differential Geom. 94, no. 1, 159--186 (2013)

\bibitem{Eichmair-Metzger2}
M.~Eichmair and J.~Metzger, \textit{On large volume preserving stable CMC surfaces in initial data sets,} J. Differential Geom. 91, no. 1, 81--102 (2012)

\bibitem{Eichmair-Metzger3}
M.~Eichmair and J.~Metzger, \textit{Unique isoperimetric foliations of asymptotically flat manifolds in all dimensions,} to appear in Invent. Math.

\bibitem{Huisken}
G.~Huisken, \textit{An isoperimetric concept for mass and quasilocal mass,} Oberwolfach reports vol. 3 no. 1, 87--88 (2012)

\bibitem{Huisken-Yau}
G.~Huisken and S.-T.~Yau, \textit{Definition of center of mass for isolated physical systems and unique foliations by stable spheres with constant mean curvature,} Invent. Math. 124, 281--311 (1996)

\bibitem{Lopez-Ros} J.~L\'opez and A.~Ros, \textit{ Complete minimal surfaces with index one and stable constant mean curvature surfaces,} Comment. Math. Helv. 64 (1989)

\bibitem{Pacard-Xu}
F.~Pacard and X.~Xu, \textit{Constant mean curvature spheres in Riemannian manifolds,} Manuscripta Math. 128, 275--295 (2008)

\bibitem{Palmer}
B.~Palmer, \textit{Surfaces of constant mean curvature in space forms,} PhD thesis, Stanford University (1986)

\bibitem{Qing-Tian}
J.~Qing and G.~Tian, \textit{On the uniqueness of the foliation of spheres of constant mean curvature in asymptotically flat $3$-manifolds,} J. Amer. Math. Soc. 20, 1091--1110 (2007)
\end{thebibliography}
\end{document}